\documentclass[11pt,a4paper]{article}
\usepackage{epsfig}
\usepackage{amsthm}

\usepackage[margin=1.1in]{geometry}
\usepackage{comment}
\usepackage{subfig}

\usepackage{amssymb,amsmath,latexsym,enumerate,graphicx,microtype,cite}
\usepackage{abstract}
\usepackage{url}
\usepackage{hyperref}
\allowdisplaybreaks

\usepackage{amsmath,amssymb,amsthm,mathrsfs}
\usepackage{mathtools}
\usepackage{comment}
\usepackage{abstract}

\usepackage[colorinlistoftodos]{todonotes}

\renewcommand{\leq}{\leqslant}
\renewcommand{\geq}{\geqslant}

\newcommand{\F}{{\cal{F}}} 
\newcommand{\A}{{\cal{A}}}

\newcommand{\bbox}{\vrule height7pt width4pt depth1pt}

\newtheorem{theorem}{Theorem}

\newtheorem{conjecture}{Conjecture}
\newtheorem{lemma}{Lemma}

\newtheorem{observation}{Observation}

\author{
	Eyal Ackerman\thanks{Department of Mathematics, Physics and Computer Science,
		University of Haifa at Oranim, 	Tivon 36006, Israel. \texttt{ackerman@math.haifa.ac.il}}\and
	G\'abor Dam\'asdi\thanks{HUN-REN Alfréd Rényi Institute of Mathematics and ELTE Eötvös Loránd University, Budapest, Hungary. Partially supported by ERC Advanced Grant GeoScape.
\texttt{damasdigabor@caesar.elte.hu}}\and
	Bal\'azs Keszegh\thanks{HUN-REN Alfréd Rényi Institute of Mathematics and ELTE Eötvös Loránd University, Budapest, Hungary. 
Research supported by the J\'anos Bolyai Research Scholarship of the Hungarian Academy of Sciences, by the National Research, Development and Innovation Office -- NKFIH under the grant K 132696 and FK 132060, by the \'UNKP-23-5 New National Excellence Program of the Ministry for Innovation and Technology from the source of the National Research, Development and Innovation Fund and by the ERC Advanced Grant ``ERMiD''. This research has been implemented with the support provided by the Ministry of Innovation and Technology of Hungary from the National Research, Development and Innovation Fund, financed under the  ELTE TKP 2021-NKTA-62 funding scheme. \texttt{keszegh@renyi.hu}}\and
	Rom Pinchasi\thanks{Technion - Israel Institute of Technology, Haifa, Israel. Visiting professor at EPFL, Lausanne, Switzerland. Supported by ISF grant (grant No.\ 1091/21). \texttt{room@technion.ac.il}}\and
	Rebeka Raffay\thanks{\texttt{rebeka.raffay@epfl.ch}}
}

\title{On the number of digons in arrangements of pairwise intersecting circles} 

\begin{document}
	\maketitle
	
	\begin{abstract}
	A long-standing open conjecture of Branko Grünbaum from 1972 states that any simple arrangement of $n$ pairwise intersecting pseudocircles in the plane can have at most $2n-2$ digons.  Agarwal et al. proved this conjecture for arrangements of pairwise intersecting pseudocircles in which there is a common point surrounded by all pseudocircles. Recently, Felsner, Roch and Scheucher showed that Gr\"unbaum's conjecture is true for arrangements of pairwise intersecting pseudocircles in which there are three pseudocircles every pair of which create a digon. In this paper we prove this over 50-year-old conjecture of Gr\"unbaum for any simple arrangement of pairwise intersecting circles in the plane.  
		
	\end{abstract}

	\section{Introduction}

	A family of \emph{pseudocircles} is a set of closed Jordan curves such that every two of them are either disjoint, intersect at exactly one point in which they touch or intersect at exactly two points in which they properly cross each other.
The bounded regions whose boundaries are the pseudocircles are called \emph{pseudodiscs}.
An \emph{arrangement} $\A(\F)$ of a family $\F$ of pseudocircles is the cell complex into which the plane is decomposed by the pseudocircles and consists of \emph{vertices}, \emph{edges} and \emph{faces}.
If there are two points that lie on every pseudocircle, then the arrangement is \emph{trivial}.
If there is no point that lies on three pseudocircles, then the arrangement is \emph{simple}.

A \emph{digon} is a face in $\A(\F)$ whose boundary consists of two edges. The two circles containing the two edges of a digon are said to
\emph{support} the digon. We also say that these two circles \emph{create} the digon.

\begin{figure}[ht]
	\centering
	\includegraphics[height=4cm]{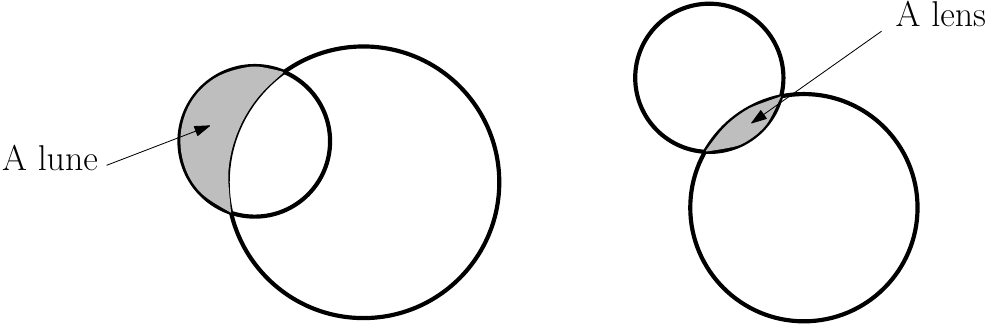}
	\caption{A lens and a lune.}
	\label{fig:lens_lune}
\end{figure}

We distinguish two different types of digons. A \emph{lens}
is a digon that is equal to the intersection of the two discs
supporting it. A \emph{lune} is a digon that is equal to a difference of the two discs supporting it (see Figure \ref{fig:lens_lune}).

It is easy to see that there are $2n$ digons in a trivial arrangement of $n$ pseudocircles, for $n>1$.
More than 50 years ago Gr\"unbaum conjectured that non-trivial arrangements of \emph{pairwise intersecting} pseudocircles have fewer digons.

\begin{conjecture}[Gr\"unbaum’s digon conjecture~{\cite[Conjecture 3.6]{GRUNB72}}]\label{conj:grunb}
	Every simple arrangement of $n > 2$ pairwise intersecting pseudocircles has at most $2n-2$ digons.
\end{conjecture}

It is possible to show, by small perturbation the pseudocircles near intersection points of three or more curves, that one can assume the family of pseudocircles is simple without decreasing the number of digons in the arrangement, as long as we do not start with a trivial arrangement.
By assuming that the arrangement of pseudocircles is 
simple, we conclude that it is nontrivial for $n>2$.

Some special cases of Gr\"unbaum's conjecture were settled.
Agarwal et al.~\cite{ANPPSS04} proved the conjecture for \emph{cylindrical} arrangements, that is, for arrangements in which there is a point that is surrounded by every pseudocircle.
Recently, Felsner, Roch and Scheucher~\cite{FELSNER23} showed that the conjecture also holds for simple arrangements in which there are three pseudocircles such that every two of them create a digon.

In this paper we prove Gr\"ubnaum's conjecture for any simple arrangement of pairwise intersecting circles in the plane.

\begin{theorem}\label{theorem:main}
	Every non-trivial simple arrangement of $n$ pairwise intersecting circles has at most $2n-2$ digons.
\end{theorem}

\begin{figure}[ht]
	\centering
	\includegraphics[height=6cm]{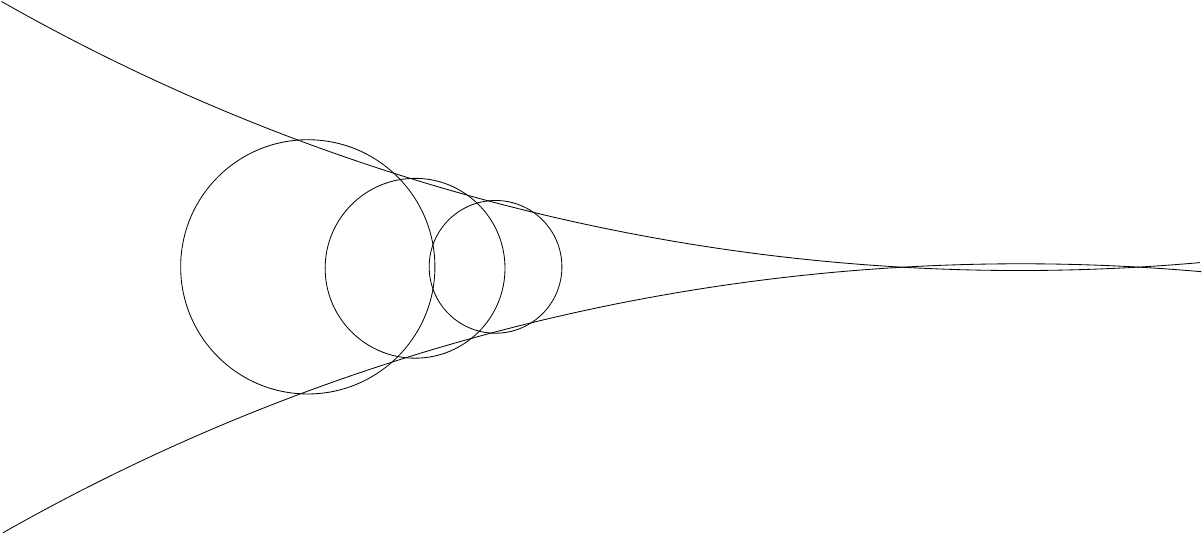}
	\caption{A family of $5$ pairwise intersecting circles with $8$ lenses.}
	\label{fig:tight}
\end{figure}

The simple construction in Figure \ref{fig:tight} (taken from~\cite{GRUNB72}) shows 
that the bound in Theorem \ref{theorem:main} is best possible for $n \geq 4$. There are $5$ circles in this construction and $8$ lenses. 
One can generalize the construction for any number of circles by suitably adding more circles to the three smaller circles
in the figure.

Before we continue we would like to note that with more care one could prove Theorem 
\ref{theorem:main} only under the assumption that the family of circles is nontrivial, rather than simple. Here, unlike in the pseudocircles case, one can no longer perturb the circles near multiple intersection points and hope to remain with a family of circles. Therefore, it will be easier and more natural to address the case of non-simple arrangements as part of the more general conjecture of Gr\"unbaum for non-simple arrangements of pairwise intersecting pseudocircles. We decided not to address in this paper the rather technical issues that may arise if one does not assume simplicity of the arrangement of circles. 
In this paper we would like to emphasize the 
special and beautiful properties of geometric circles related to digons. 
The assumption that the family of circles is simple allows
one to assume, and we will indeed assume this, that the set of centers of the circles in question is in general position in the sense that no three of the centers are collinear.
This can be done by applying a generic inversion map to the plane. We will not make use these assumptions explicitly in the proofs, but we will use them implicitly.

There has been a lot of research about digons in arrangements of circles (and pseudo-circles) that are not necessarily pairwise intersecting. We will not survey here the vast literature about digons in arrangements of circles and pseudo-circles and on related situations where we allow curves to intersect more than twice and only refer the reader to \cite{GOT18} and the many references therein. The case where circles need not be pairwise intersecting is of completely different nature. We remark that in such a case the best constructions show that it is possible that $n$ circles will determine $\Omega(n^{4/3})$ many lenses. The best known upper bound is $O(n^{3/2}\log n)$ given in \cite{MT06}, that is following the footsteps of \cite{PR04}. A slightly better upper bound of $O(n^{3/2})$ for the number of touching points among $n$ circles follows from a result of Ellenberg, Solymosi and Zahl~\cite{ESZ16}. The case of unit circles is of particular interest because of its relation to the celebrated unit distance problem posed by Paul Erd\H os (\cite{Erdos46}). For this problem the best known lower and upper bounds are $\Omega(n^{1+c/\log\log n})$~\cite{Erdos46} and $O(n^{4/3})$~\cite{PT06,SST84,Szekely97}, respectively.   

Going back to families of pairwise intersecting circles, the number of lunes in these arrangements was studied in \cite{ALPS01}. 

\begin{theorem}\label{theorem:lunes}
Any arrangement of $n$ pairwise intersecting circles in the plane 
has at most $2n-4$ lunes.
\end{theorem}

Theorem \ref{theorem:lunes} is used in \cite{ALPS01}
to derive a linear upper bound, that is not tight, for the number of digons (lunes and lenses) in any arrangement of pairwise intersecting circles in the plane. 
Specifically, it is shown in \cite{ALPS01} that arrangements of $n$ pairwise intersecting circles in the plane contain at most $2n-2$ lunes and at most $18n$ lenses.

The tight bound on the maximum number of lenses in a family of pairwise intersecting circles in the plane is established in \cite{P24}.

\begin{theorem}[\cite{P24}]\label{theorem:lenses}
Any arrangement of $n$ pairwise intersecting circles in the plane determines at most $2n-2$ lenses.
\end{theorem}

Theorem \ref{theorem:lunes} and Theorem \ref{theorem:lenses} imply immediately an upper bound of $4n-6$ for the number of digons 
in arrangements of $n$ pairwise intersecting
circles in the plane.
Hence, Theorem~\ref{theorem:main} improves this upper bound.

For arrangements of pairwise intersecting \emph{unit circles}, Pinchasi~\cite{Pinchasi02} proved that they can have at most $n$ lenses and at most $3$ lunes, hence at most $n+3$ digons.

\section{Three crucial geometric lemmata}

In this section we will explore three crucial geometric lemmata concerning pairwise intersecting circles in the plane. Two of these lammata have been shown in previous works
while the third one is new and we bring its proof here.

All the three lemmata are concerned with 
the geometric graph $G$ on the set of centers of circles in a family $\F$ of pairwise intersecting circles in the plane. The edges in $G$ correspond to the pairs of circles creating digons (either lunes, or lenses, or both) in $\A(\F)$.

The first lemma is from \cite{ALPS01}, where it is used to derive the tight upper bound for the number of lunes in Theorem \ref{theorem:lunes}.

\begin{lemma}[\cite{ALPS01}]\label{lemma:lunes}
Let $\F$ be a finite family of pairwise intersecting circles in the plane. Let
$G$ be the geometric graph whose vertices are the centers of the circles in $\F$
such that two vertices are connected by an edge if and only if the corresponding circles in $\F$
create a lune in the arrangement $\A(\F)$.
Then no pair of edges in $G$ cross, thus $G$ is a planar embedding.
\end{lemma}

\begin{figure}[ht]
	\centering
	\includegraphics[height=3cm]{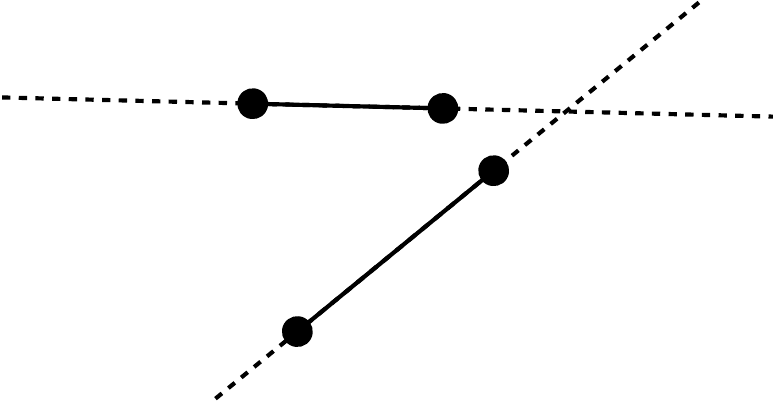}
	\caption{A pair of avoiding edges.}
	\label{fig:avoiding}
\end{figure}

The second lemma is from \cite{P24},
where it is used to prove Theorem \ref{theorem:lenses}.
We recall that two edges in a geometric graph are called \emph{avoiding} if they are opposite edges of a convex quadrilateral. That is, two edges $e$ and $f$ in a geometric graph are avoiding if the lines containing $e$ and $f$ intersect at a point that is neither on $e$ nor nor on $f$ (see Figure \ref{fig:avoiding}).

\begin{lemma}[\cite{P24}]\label{lemma:lenses}
Let $\F$ be a finite family of pairwise intersecting circles in the plane. Let
$G$ be the geometric graph whose vertices are the centers of the circles in $\F$
such that two vertices are connected by an edge if and only if the corresponding circles in $\F$
create a lens in the arrangement $\A(\F)$.
Then $G$ does not contain a pair of avoiding edges.
\end{lemma}

Lemma \ref{lemma:lenses} implies immediately Theorem \ref{theorem:lenses}
because of a result of Katchalski and Last \cite{KL98} and Valtr \cite{V98} by which a geometric graph on $n$ vertices with no pair of avoiding edges
can have at most $2n-2$ edges. We will not make use of this bound, but rather use the observation in Lemma \ref{lemma:lenses} directly.

For the sake of completeness and because 
the result in \cite{P24} is rather recent, we bring here an independent proof, and different from the one in \cite{P24}, of Lemma \ref{lemma:lenses} in the case where 
the arrangement of pairwise intersecting circles is simple. 

\bigskip

\noindent {\bf Proof of Lemma
\ref{lemma:lenses}.}
Under the contrary assumption, there are four circles $c_{1}, c_{2}, c_{3}$, and $c_{4}$ in $\F$ 
with centers $O_{1}, O_{2}, O_{3}$, and $O_{4}$, respectively,
such that $c_{1}$ and $c_{2}$ create a lens in $\A(\F)$ and also $c_{3}$ and $c_{4}$ create a lens in $\A(\F)$. Moreover, the line segments 
$[O_{1}O_{2}]$ and $[O_{3}O_{4}]$ are 
opposite edges of a convex quadrilateral.

\begin{figure}[ht]
	\centering
	\includegraphics[height=6cm]{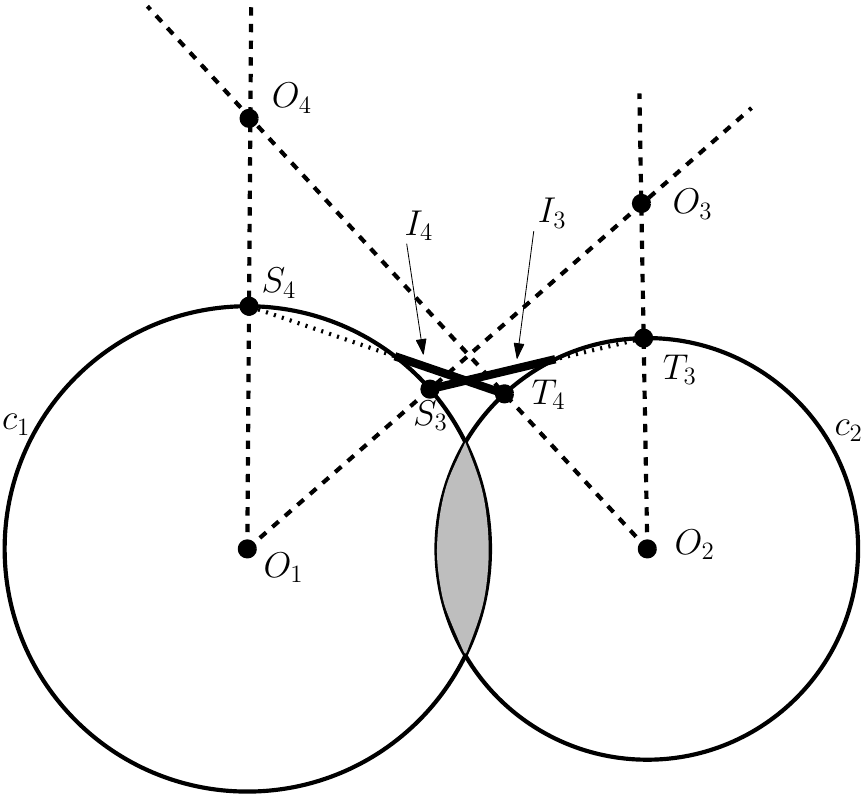}
	\caption{Lemma \ref{lemma:lenses}.}
	\label{fig:lemma_lenses}
\end{figure}

Without loss of generality we assume that the line $O_{1}O_{2}$ is horizontal and
$O_{1}$ lies to the left of $O_{2}$. We may also assume that both $O_{3}$ and $O_{4}$ lie 
above the line $O_{1}O_{2}$ such that 
$O_{1}O_{2}O_{3}O_{4}$ is a convex quadrilateral (see Figure \ref{fig:lemma_lenses}).

We need the following two simple observations.
The first is extremely elementary and its proof is left to the reader. The second observation is a common knowledge that is very well known.

\begin{figure}[ht]
	\centering
	\includegraphics[height=3.5cm]{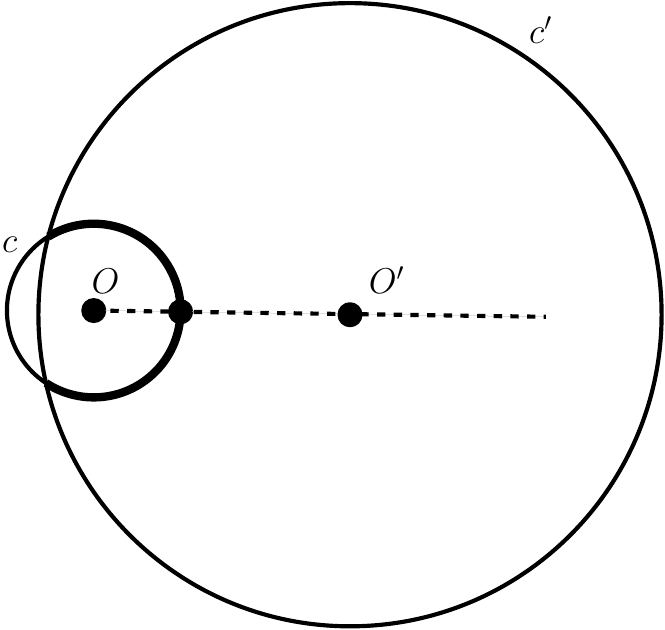}
	\caption{An illustration for Observation \ref{observation:simple}.}
	\label{fig:simple}
\end{figure}

\begin{observation}\label{observation:simple}
Let $c$ and $c'$ be two intersecting circles
with centers $O$ and $O'$ respectively.
Then the intersection of $c$ and the ray $\overrightarrow{OO'}$ is the center of the arc
on $c$ that is the part of $c$ inside the disc bounded by $c'$ (see Figure \ref{fig:simple}).
\end{observation}

\begin{observation}\label{observation:d}
Let $A,B,A'$, and $B'$ be four points in the plane.
If there is a circular disc $D$ that contains $A$ and $B$ but not $A'$ and $B'$, and there is another circular disc $D'$ that contains $A'$ and $B'$ but not $A$ and $B$, then the two segments $[AB]$ and $[A'B']$ are disjoint.
\end{observation}

\begin{figure}[ht]
	\centering
	\includegraphics[height=4cm]{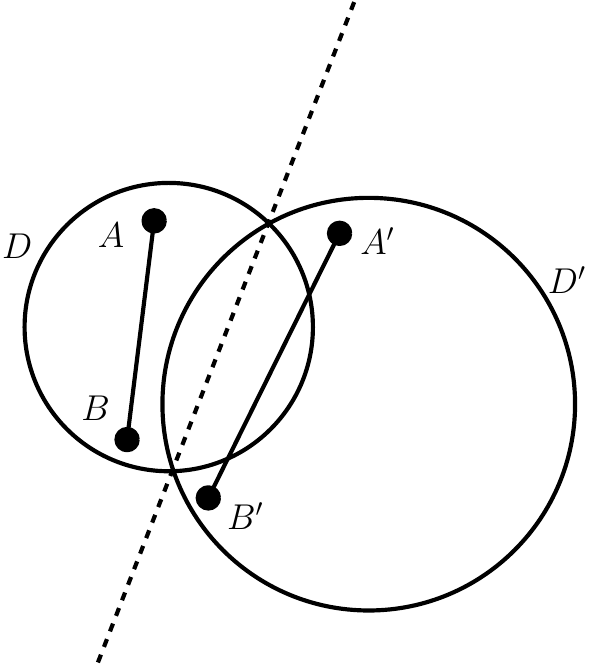}
	\caption{An illustration for Observation \ref{observation:d}.}
	\label{fig:d}
\end{figure}

\noindent {\bf Proof.} This is clear if the circular discs
$D$ and $D'$ are disjoint. If $D$ and $D'$ intersect
we observe that the line through the intersection points
of their boundaries (notice none of $D$ and $D'$ can be contained in the other)  separates the regions $D \setminus D'$ and $D' \setminus D$
(see Figure \ref{fig:d}).
\bbox

\bigskip

Going back to the proof of Lemma \ref{lemma:lenses}, let $D_{3}$ and $D_{4}$ denote the circular discs bounded by $c_{3}$ and $c_{4}$, respectively.
let $S_{3}$ and $S_{4}$ denote the centers
of the arcs $c_{1} \cap D_{3}$ and $c_{1} \cap D_{4}$, respectively, on the circle $c_{1}$. We observe that the two arcs 
must be disjoint, or else the lens $D_{3} \cap D_{4}$ would have contained a point of $c_{1}$, which is impossible.

Similarly, let $T_{3}$ and $T_{4}$
denote the centers of the arcs $c_{2} \cap D_{3}$ and $c_{2} \cap D_{4}$, respectively, on the circle $c_{2}$
and observe that the two arcs must be disjoint.

By observation \ref{observation:simple},
$S_{3}$ is the point of intersection of $c_{1}$ with $\overrightarrow{O_{2}O_{3}}$. Similarly, $S_{4}$ is the point of intersection of $c_{1}$ with $\overrightarrow{O_{2}O_{4}}$.
Because $O_{1}O_{2}O_{3}O_{4}$ is a convex quadrilateral, it must be that $S_{3}$ lies to the right of $S_{4}$ on $c_{1}$
above the line $O_{1}O_{2}$ (see Figure \ref{fig:lemma_lenses}).

We claim that $D_{4}$ cannot contain the lens created by $c_{1}$ and $c_{2}$. Indeed, $c_{1} \cap D_{4}$ is an arc whose center $S_{4}$ lies above the line $O_{1}O_{2}$. Therefore, if $D_{4}$ contains the lens created by $c_{1}$ and $c_{2}$, then $c_{1} \cap D_{4}$ must contain all the part of $c_{1}$ that is above the line $O_{1}O_{2}$ and to the right of $S_{4}$. In particular it must contain the center $S_{3}$ of the arc $c_{1} \cap D_{3}$.
Consequently, the point $S_{3}$ lies both in the interior of $D_{4}$ and the interior of $D_{3}$. This is a contradiction because the lens $D_{3} \cap D_{4}$ cannot contain in its interior a point of $c_{1}$.

We argue similarly about $T_{3}$ and $T_{4}$. By Observation \ref{observation:simple},
$T_{3}$ is the point of intersection of $c_{2}$ with $\overrightarrow{O_{2}O_{3}}$. Similarly, $T_{4}$ is the point of intersection of $c_{2}$ with $\overrightarrow{O_{2}O_{4}}$.
Because $O_{1}O_{2}O_{3}O_{4}$ is a convex quadrilateral $T_{3}$ lies to the right of $T_{4}$ on $c_{2}$
above the line $O_{1}O_{2}$.
We observe now that $D_{3}$ cannot contain the lens created by $c_{1}$ and $c_{2}$. This is because 
otherwise $D_{3}$ contains the arc $c_{2} \cap D_{4}$ whose center is the point $T_{4}$ on $c_{2}$. Hence $T_{4}$, that is a point on $c_{2}$, lies in the interiors of both $D_{4}$ and $D_{3}$. Consequently $c_{2}$ intersects the interior of the lens $D_{3} \cap D_{4}$, which is a contradiction.

Because $D_{3}$ and $D_{4}$ do not contain the lens $D_{1} \cap D_{2}$, they must be disjoint from it, as the lens $D_{1} \cap D_{2}$ cannot be intersected by any of the circles $c_{1}$ and $c_{2}$. Hence all four arcs $c_{1} \cap D_{3}, c_{1} \cap D_{4}, c_{2} \cap D_{3}$, and $c_{2} \cap D_{4}$
lie on the boundary of
$D_{1} \cup D_{2}$. This boundary is a simple closed curve that we denote by $\Gamma$. The curve $\Gamma$ is the union of the two arcs $c_{1} \setminus D_{2}$ and $c_{2} \setminus D_{1}$. 

Because $S_{3}$ lies to the right of $S_{4}$ on $c_{1}$
and $T_{3}$ lies to the right of $T_{4}$ on $c_{2}$, above the line $O_{1}O_{2}$, then 
the four pairwise disjoint arcs
$c_{2} \cap D_{3}, c_{2} \cap D_{4}, c_{1} \cap D_{3}$, and $c_{1} \cap D_{4}$ lie in this  counterclockwise cyclic order on $\Gamma$.
In particular, the arcs $c_{1} \cap D_{3}$ and $c_{2} \cap D_{3}$ separate the arcs
$c_{1} \cap D_{4}$ and $c_{2} \cap D_{4}$ on the simple closed curve $\Gamma$.

We now claim that the line segments $[S_{3}T_{3}]$ and $[S_{4}T_{4}]$
must cross. This will lead to a contradiction 
because these two segments must be disjoint
by Observation \ref{observation:d}. This is because both points
$S_{3}$ and $T_{3}$ belong to $D_{3}$ and not to $D_{4}$ while both points $S_{4}$ and $T_{4}$
belong to $D_{4}$ and not to $D_{3}$.

To see that the line segments $[S_{3}T_{3}]$ and $[S_{4}T_{4}]$
must cross, we observe that the segment 
$[S_{3}T_{3}]$ intersects the interior of the region 
$\mathbb{R}^2 \setminus D_{1} \cup D_{2}$ at a chord $I_{3}$
connecting a point on the arc $c_{1} \cap D_{3}$ (this point could be $S_{3}$, but 
not necessarily) with a point on the arc $c_{2} \cap D_{3}$. 
Similarly, $[S_{4}T_{4}]$ intersects 
the interior of  
$\mathbb{R}^2 \setminus D_{1} \cup D_{2}$ at a chord $I_{4}$
connecting a point on the arc $c_{1} \cap D_{4}$ with a point on the arc
$c_{2} \cap D_{4}$. 

It now follows that the two endpoints of the chord $I_{4}$
separate the two endpoints of the chord $I_{3}$ on the simple closed curve $\Gamma$, that is also the boundary of $\mathbb{R}^2 \setminus D_{1} \cup D_{2}$ (see Figure \ref{fig:lemma_lenses}).
Because both $I_{3}$ and $I_{4}$ are contained in $\mathbb{R}^2 \setminus D_{1} \cup D_{2}$ whose boundary is the simple closed curve $\Gamma$, 
it follows that $I_{3}$ and $I_{4}$
must cross inside the region $\mathbb{R}^2 \setminus D_{1} \cup D_{2}$. Consequently, $[S_{3}T_{3}]$ and $[S_{4}T_{4}]$ cross, which is the desired contradiction.
\bbox

\bigskip

The third and last lemma that we will need for the proof of Theorem \ref{theorem:main}
is new. Similar to 
Lemma \ref{lemma:lunes} and Lemma \ref{lemma:lenses}, this lemma is concerned too with the geometric graph representing the lunes and lenses in arrangements of pairwise intersecting circles and it is concerned with the mutual relations between lunes and lenses in such arrangements. 

\begin{figure}[ht]
	\centering
	\includegraphics[height=4cm]{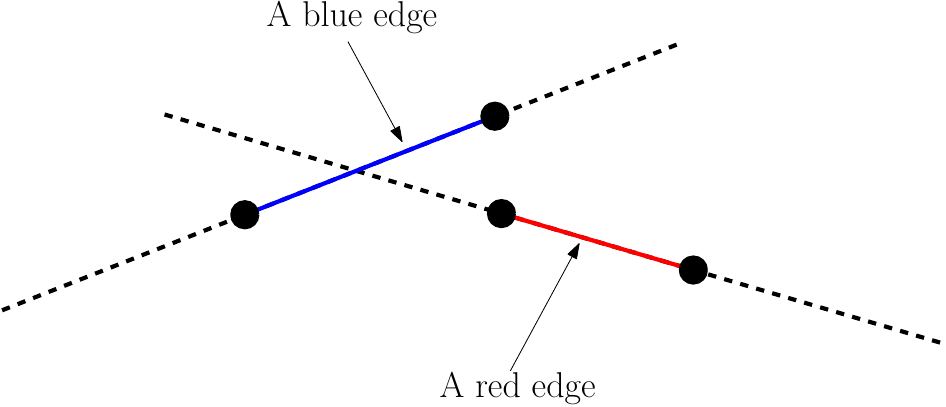}
	\caption{The forbidden position of red and blue edges.}
	\label{fig:red_blue_edges}
\end{figure}

\begin{lemma}\label{lemma:lunes_lenses}
Let $\F$ be a family of pairwise intersecting circles in the plane. Define
a geometric graph $G$ on the set of centers
of the circles in $\F$. We connect two vertices (centers) in $G$ with a blue edge
if the corresponding circles create a lune in $\A(\F)$. We connect two vertices in $G$ with a red edge if the corresponding circles create a lens in $\A(\F)$.
Then $G$ does not contain a red edge $e$ and a blue edge $f$ such that $e$ and $f$ are disjoint and the line through $e$ intersects $f$ (see Figure \ref{fig:red_blue_edges}).
\end{lemma}

\noindent{\bf Proof.}
Under the contrary assumption, there are four circles $c_{1}, c_{2}, c_{3}$, and $c_{4}$ in $\F$ 
with centers $O_{1}, O_{2}, O_{3}$, and $O_{4}$, respectively,
such that $c_{1}$ and $c_{2}$ create a lune in $\A(\F)$ while $c_{3}$ and $c_{4}$ create a lens in $\A(\F)$. Moreover, the line segments 
$[O_{1}O_{2}]$ and $[O_{3}O_{4}]$ are disjoint and the line $O_{3}O_{4}$ intersects 
the line segment $[O_{1}O_{2}]$.

\begin{figure}[ht]
	\centering
	\includegraphics[height=6cm]{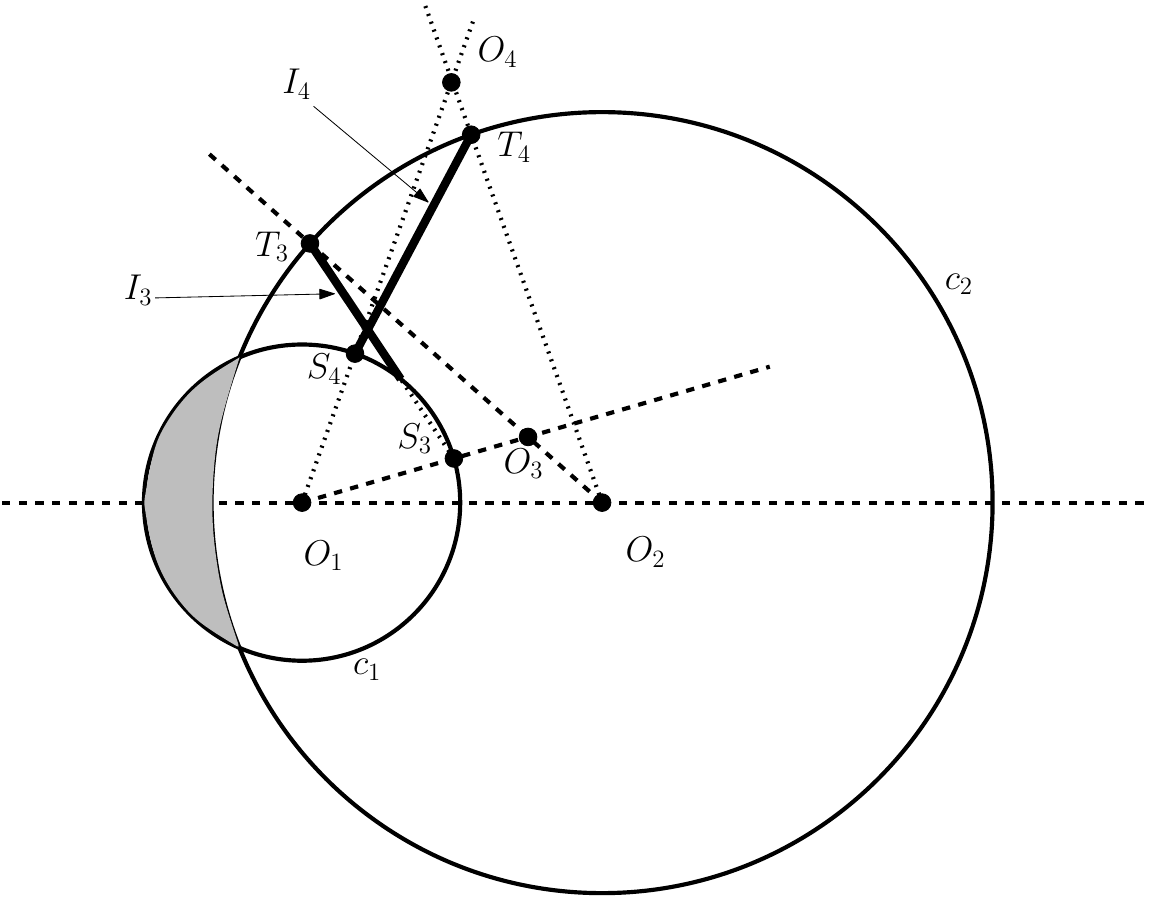}
	\caption{Lemma \ref{lemma:lunes_lenses}.}
	\label{fig:lemma_lens_lune}
\end{figure}

Without loss of generality we assume that 
the lune created by $c_{1}$ and $c_{2}$ is equal to the disc bounded by $c_{1}$ minus the disc bounded by $c_{2}$. We can also assume without loss of generality that $O_{1}O_{2}$ is horizontal and
$O_{1}$ lies to the left of $O_{2}$. We may also assume that both $O_{3}$ and $O_{4}$ lie 
above the line $O_{1}O_{2}$ and $O_{3}$ is closer than $O_{4}$ to the line $O_{1}O_{2}$. By our assumptions, the point $O_{4}$ must belong to the angle opposite to $\angle O_{1}O_{3}O_{2}$ (see Figure \ref{fig:lemma_lens_lune}).

Let $D_{3}$ and $D_{4}$ denote the circular discs bounded by $c_{3}$ and $c_{4}$, respectively.
Let $S_{3}$ and $S_{4}$ denote the centers
of the arcs $c_{1} \cap D_{3}$ and $c_{1} \cap D_{4}$, respectively, on the circle $c_{1}$. We observe that the two arcs 
must be disjoint, or else the lens $D_{3} \cap D_{4}$ would have contained a point of $c_{1}$.
Similarly, let $T_{3}$ and $T_{4}$
denote the centers of the arcs $c_{2} \cap D_{3}$ and $c_{2} \cap D_{4}$, respectively, on the circle $c_{2}$
and observe that the two arcs must be disjoint.

By observation \ref{observation:simple},
$S_{3}$ is the point of intersection of $c_{1}$ with $\overrightarrow{O_{2}O_{3}}$. Similarly, $S_{4}$ is the point of intersection of $c_{1}$ with $\overrightarrow{O_{2}O_{4}}$.
Because $O_{4}$ belongs to the angle opposite to $\angle O_{1}O_{3}O_{2}$, it must be that $S_{3}$ lies to the right of $S_{4}$ on $c_{1}$
above the line $O_{1}O_{2}$.

We observe now that $D_{3}$ cannot contain the lune created by $c_{1}$ and $c_{2}$. 
Indeed, $c_{1} \cap D_{3}$ is an arc whose center $S_{3}$ lies above the line $O_{1}O_{2}$. Therefore, if $D_{3}$ contains the lune created by $c_{1}$ and $c_{2}$, then $c_{1} \cap D_{3}$ must contain all the part of $c_{1}$ that is above the line $O_{1}O_{2}$ and to the left of $S_{3}$. In particular it contains the arc $c_{1} \cap D_{4}$ and its center $S_{4}$.
Then $S_{4}$ is contained in the interiors of both $D_{3}$ and $D_{4}$ and therefore it is contained in the interior of the lens $D_{3} \cap D_{4}$. This is impossible because $S_{4}$ is a point on $c_{1}$ that must be disjoint from the lens $D_{3} \cap D_{4}$.

We now argue similarly about $T_{3}$ and $T_{4}$. By Observation \ref{observation:simple},
$T_{3}$ is the point of intersection of $c_{2}$ with $\overrightarrow{O_{2}O_{3}}$. Similarly, $T_{4}$ is the point of intersection of $c_{2}$ with $\overrightarrow{O_{2}O_{4}}$.
Because $O_{4}$ belongs to the angle opposite to $\angle O_{1}O_{3}O_{2}$, it must be that $T_{4}$ lies to the right of $T_{3}$ on $c_{2}$
above the line $O_{1}O_{2}$.
We observe now that $D_{4}$ cannot contain the lune created by $c_{1}$ and $c_{2}$. This is because 
otherwise $D_{4}$ contains the point $T_{3}$ on $c_{2}$.
Consequently $T_{3}$ is contained in the interiors of both $D_{4}$ and $D_{3}$. This is impossible because $T_{3}$ is a point on $c_{2}$ that is disjoint from the lens $D_{3} \cap D_{4}$.

Because $D_{3}$ and $D_{4}$ do not contain the lune $D_{1} \setminus D_{2}$, they must be disjoint from it, as a lune $D_{1} \setminus D_{2}$ cannot be intersected by any of the circles $c_{1}$ and $c_{2}$. Hence all four arcs $c_{1} \cap D_{3}, c_{1} \cap D_{4}, c_{2} \cap D_{3}$, and $c_{2} \cap D_{4}$
lie on the boundary of
$D_{2} \setminus D_{1}$. This boundary is a simple closed curve that we denote by $\Gamma$.

It follows now that the pairwise disjoint arcs $c_{2} \cap D_{4}, c_{2} \cap D_{3}, c_{1} \cap D_{4}$, and $c_{1} \cap D_{3}$ lie in this counterclockwise cyclic order on the closed curve $\Gamma$.
Notice in particular that the arcs $c_{2} \cap D_{3}$ and $c_{1} \cap D_{3}$ separate the arcs $c_{2} \cap D_{4}$ and $c_{1} \cap D_{4}$ on the simple closed curve $\Gamma$.

We claim that the line segments $[S_{3}T_{3}]$ and $[S_{4}T_{4}]$
must cross. This will lead to a contradiction 
because these two segments must be disjoint
by Observation \ref{observation:d}. This is because both points
$S_{3}$ and $T_{3}$ belong to $D_{3}$ and not to $D_{4}$ while both points $S_{4}$ and $T_{4}$
belong to $D_{4}$ and not to $D_{3}$.

To see that the line segments $[S_{3}T_{3}]$ and $[S_{4}T_{4}]$
must cross, we observe that the segment 
$[S_{3}T_{3}]$ intersects the interior of the simply connected region 
$D_{2} \setminus D_{1}$ at a chord $I_{3}$
connecting a point on the arc $c_{1} \cap D_{3}$ (this point could be $S_{3}$, but 
not necessarily) with the point $T_{3}$ on the arc $c_{2} \cap D_{3}$. 

Similarly, $[S_{4}T_{4}]$ intersects 
the interior of the region 
$D_{2} \setminus D_{1}$ at a chord $I_{4}$
connecting a point on the arc $c_{1} \cap D_{4}$ with the point $T_{4}$ on the arc
$c_{2} \cap D_{4}$. In particular, the two endpoints of the chord $I_{4}$
separate the two endpoints of the chord $I_{3}$ on the simple closed curve $\Gamma$.
Because both $I_{3}$ and $I_{4}$ are contained in simply connected region $D_{2} \setminus D_{1}$ whose boundary is $\Gamma$,
then $I_{3}$ and $I_{4}$
must cross inside the region $D_{2} \setminus
D_{1}$ (see Figure \ref{fig:lemma_lens_lune}). Consequently, $[S_{3}T_{3}]$ and $[S_{4}T_{4}]$ cross, which is the desired contradiction.
\bbox

\bigskip

The proof of Lemma \ref{lemma:lunes_lenses} resembles a lot the proof of Lemma \ref{lemma:lenses} that we brought above.
We expect that one can state and prove a unified version of all three lemmata presented in this section for arrangements of pairwise intersecting circles on the sphere. 

\section{Proof of Theorem \ref{theorem:main}}

We consider the geometric graph $G$ as in the statement of Lemma 
\ref{lemma:lunes_lenses}. That is, the vertices of $G$ are the centers of the circles in $\F$. Two centers are connected with a red edge if the corresponding circles
create a lens. Two centers are connected with a blue edge if the corresponding circles create a lune.

We may assume without loss of generality that every circle in $\F$ supports some digon, which is either a lens or a lune.
For a circle $c$ in $\F$ we say that it is \emph{internal} if it supports a digon that is surrounded by $c$. 
We say that $c$ is \emph{external} if it supports a digon (necessarily a lune) that is not surrounded by $c$.

The following observation is simple and yet important for the proof.

\begin{observation}\label{observation:int_ext}
A circle in $\F$ cannot be both internal and external.
\end{observation}

\noindent {\bf Proof.}
Assume to the contrary that $c \in \F$ is both external and internal. Then there is a circle $c_{1} \in \F$ such that $c_{1}$ and $c$ create a digon (must be a lune) that is not surrounded by $c$. Similarly, there is a circle $c_{2} \in \F$ such that $c$ and $c_{2}$ create a digon (could be a lune or a lens) that is surrounded by $c$.

\begin{figure}[ht]
	\centering
	\includegraphics[height=4cm]{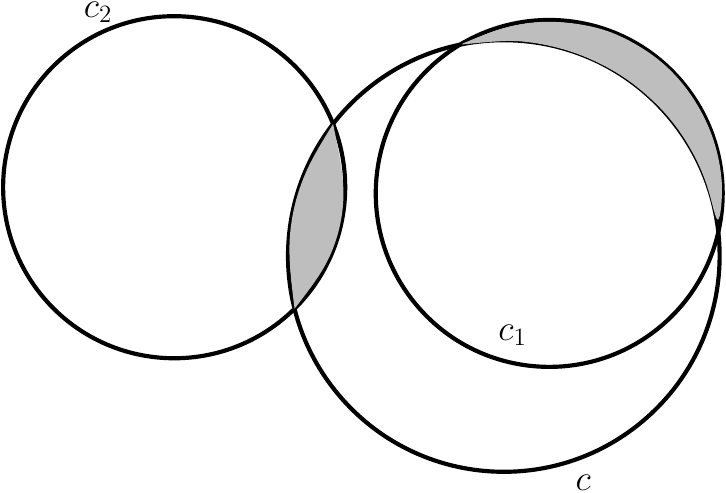}
	\caption{A circle cannot be both internal and external.}
	\label{fig:int_ext}
\end{figure}

We show that $c_{1}$ and $c_{2}$ cannot cross and thus obtain a contradiction because $\F$ is a family of pairwise intersecting circles (see Figure \ref{fig:int_ext}).
Indeed, $c_{1}$ and $c_{2}$ cannot cross inside the disc bounded by $c$ because 
the part of $c_{2}$ there is an edge of a face (digon) in $\A(\F)$. In very much the same way
$c_{1}$ and $c_{2}$ cannot cross in the region not bounded by $c$ because the part of $c_{1}$ there is an edge of a face (in fact a lune) in $\A(\F)$.
\bbox

\bigskip

As a consequence of Observation \ref{observation:int_ext}, we see that 
the blue edges connect centers of internal circles to centers of external circles, thus forming a bipartite graph. The red edges in $G$ connect centers of pairs
of internal circles.

\begin{figure}[ht]
	\centering
	\includegraphics[height=4cm]{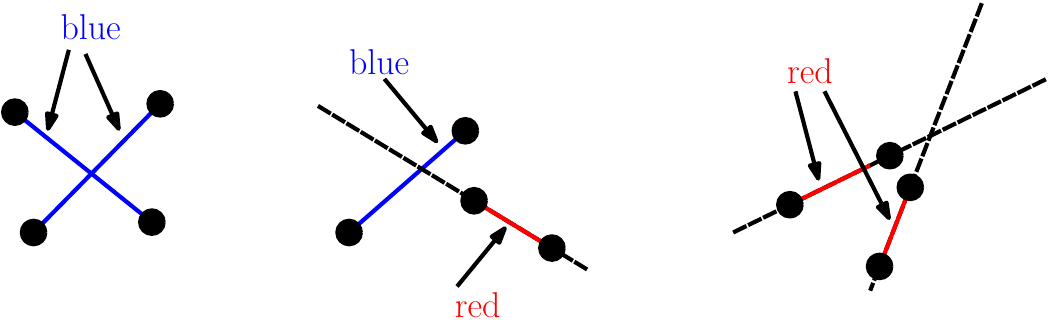}
	\caption{The forbidden pairs of edges in $G$.}
	\label{fig:forbidden}
\end{figure}

By Lemma \ref{lemma:lunes}, no two blue edges cross each other. By Lemma \ref{lemma:lenses}, no two red edges are avoiding. By Lemma \ref{lemma:lunes_lenses}, it is not possible that the complement of a red edge on the line containing it crosses a blue edge (see Figure \ref{fig:forbidden}).

With an aid of a nice trick we will move from the graph $G$ to another graph $G'$
that has twice as many edges as $G$ and twice as many vertices. The graph $G'$ will be bipartite and planar. This will show that $2E \leq 2(2n)-4$, where $E$ is the number of edges in $G$. This implies $E \leq 2n-2$ as desired.

\begin{figure}[ht]
	\centering
	\includegraphics[height=6cm]{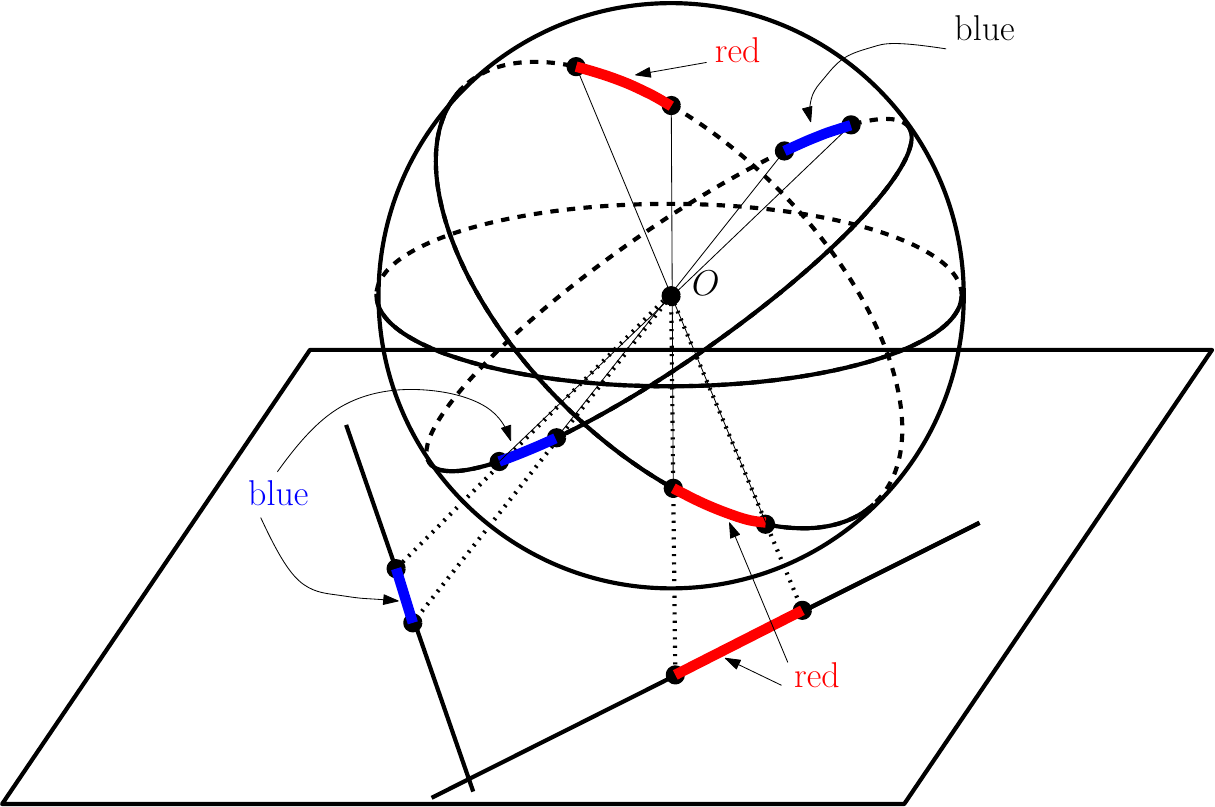}
	\caption{Projecting $G$ on the southern hemisphere and reflecting on the northern hemisphere.}
	\label{fig:trick_2}
\end{figure}

To construct $G'$ we place a sphere $S$ touching the plane from above, thinking of the plane as horizontal and the center of $S$ as the origin $O$.
Then we use a central projection from the center of $S$ and project the plane including the drawing of $G$ to the southern hemishpere of $S$ (see Figure \ref{fig:trick_2}).

We get a drawing of $G$ on the southern hemisphere of $S$ where the edges are great arcs of $S$. Next we duplicate the drawing on the northern hemisphere
by reflecting the southern hemisphere through the center of $S$. That is,
we take the drawing of $G$ on the southern hemisphere of $S$ and also take the minus of this drawing with respect to the origin $O$ that is also the center of $S$ (see Figure \ref{fig:trick_2}).

We thus get two drawings of the graph $G$ on $S$. One on the southern hemisphere and one on the northern hemisphere (its minus).

\begin{figure}[ht]
	\centering
	\includegraphics[height=7cm]{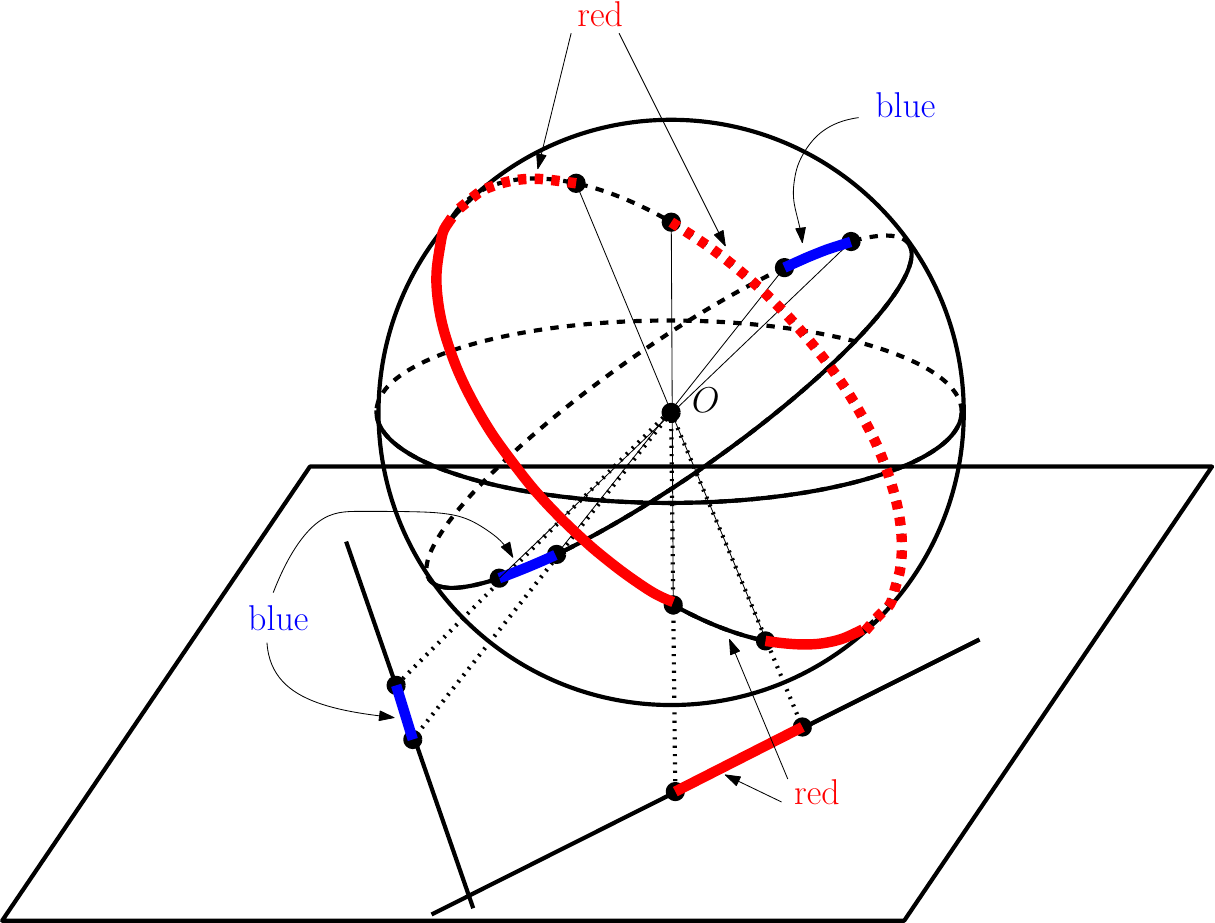}
	\caption{Replacing the red edges by their complements on the great circle.}
	\label{fig:trick_3}
\end{figure}

Next we perform the following change in the drawing to get the graph $G'$ drawn on $S$.
For every red edge $e$ in the southern hemisphere we replace $e$ and its reflection $-e$ by the complementary arcs
on the great circle on $S$ containing them (see Figure \ref{fig:trick_3}).

Therefore, $G'$ is a graph on $2n$ vertices
drawn on $S$ with twice as many red and blue edges as in the graph $G$.

The nice and crucial observation is that no two edges of $G'$ may cross.
This is directly follows from the three lemmata Lemma \ref{lemma:lunes}, Lemma \ref{lemma:lenses}, and Lemma \ref{lemma:lunes_lenses}. Moreover,
the graph $G'$ is easily seen to be bipartite. Indeed, denote by $A$ the set of external vertices on the southern hemisphere of $S$ and denote by $B$ the set of internal vertices on the southern hemisphere of $S$. Then the blue edges in $G'$ run between $A$ and $B$ and between $-A$ and $-B$. The red edges in $G'$ run
between vertices in $B$ and vertices in $-B$. Therefore, $G'$ is bipartite with the two parts being $A \cup -B$ and $B \cup -A$. This completes the proof of Theorem \ref{theorem:main}.
\bbox

\bibliographystyle{abbrv}

\begin{thebibliography}{10}
	
\footnotesize

\bibitem{ANPPSS04} 
{P. Agarwal, E. Nevo, J. Pach, R. Pinchasi, M. Sharir,
and S. Smorodinsky,}
{Lenses in Arrangements of Pseudocircles and their Applications,}
{\it J. ACM}, {\bf 51}, 139--186,  (2004) 

\bibitem{ALPS01} {N. Alon, H. Last, R. Pinchasi, and M. Sharir}, 
{On the Complexity of Arrangements of Circles in the Plane},
{\it Discrete and Computational Geometry}, {\bf 26}, 465--492, (2001)

\bibitem{ESZ16}{J. Ellenberg, J. Solymosi, and J. Zahl},  New bounds on curve tangencies and orthogonalities. {\em Discrete Anal.}, {\bf 18},  22, (2016)

\bibitem{Erdos46}
{P. Erd\H os,} 
On sets of distances of $n$ points {\it Amer. Math. Monthly.} {\bf 53}, 248--250, (1946)

\bibitem{FELSNER23}
{S. Felsner, S. Roch, M. Scheucher},
{Arrangements of Pseudocircles: On Digons and Triangles}, {\em Graph Drawing And Network Visualization}. \textbf{13764}  441-455 (2023), 


\bibitem{FELSNER21}{S. Felsner, M. Scheucher},
 Arrangements of pseudocircles: triangles and drawings. {\em Discrete Comput. Geom.}. \textbf{65}, 261-278, (2021)

\bibitem{GOT18}
Handbook of discrete and computational geometry. Third edition. Edited by Jacob E. Goodman, Joseph O'Rourke and Csaba D. T\'oth. Discrete Mathematics and its Applications (Boca Raton). CRC Press, Boca Raton, FL, (2018)

\bibitem{GRUNB72} 
{B. Gr\"unbaum},
{Arrangements and Spreads}, CBMS Regional Conference Series in
Mathematics, vol. 10. AMS (1972)

 \bibitem{KL98} {M. Katchalski, H. Last}, On geometric graphs with no two edges in convex position. Dedicated to the memory of Paul Erd\H os. {\em Discrete Comput. Geom.} {\bf 19}  no. 3, Special Issue, 399--404, (1998)
 
\bibitem{MT06} {A. Marcus, G. Tardos}, Intersection reverse sequences and geometric applications. {\it J. Combin. Theory Ser. A} {\bf 113} no. 4, 675--691, (2006)

\bibitem{PT06} {J. Pach, G. Tardos}, Forbidden paths and cycles in ordered graphs and matrices. {\em Isr. J. Math.}. \textbf{155}, 359-380, (2006)

 
\bibitem{P08} {R. Pinchasi},
{Geometric graphs with no two parallel edges},
{\it Combinatorica}, {\bf 28} no. 1, 127--130, (2008)

\bibitem{P24}{R. Pinchasi},
{A note on lenses in arrangements of pairwise intersecting circles in the plane}.
{\it Elec. J. of Comb}, accepted.

\bibitem{Pinchasi02} 
{R. Pinchasi},
{Gallai-Sylvester Theorem for Pairwise Intersecting Unit 
Circles},
{\it Discrete and Computational Geometry},
{\bf 28} 607--624, (2002)

\bibitem{PR04} 
{R. Pinchasi and R. Radoi\v ci\'c},
{On the Number of Edges in a Topological Graph with no 
Self-intersecting Cycle of Length $4$},
{\it Towards a Theory of Geometric Graphs}, 233--243, Contemp. Math., 342,
Amer. Math. Soc. Providence, RI, (2004)

\bibitem{SST84} {J. Spencer, E. Szemerédi, and W. Trotter},
Unit distances in the Euclidean plane. {\em Graph Theory And Combinatorics}. 294-304, (1984)

\bibitem{Szekely97} {L. Székely},
Crossing Numbers and Hard Erdos Problems in Discrete Geometry. {\em Comb. Probab. Comput.}. \textbf{6}, 353-358, (1997) 

\bibitem{V98} {P. Valtr}, On geometric graphs with no $k$ pairwise parallel edges, {\em Discrete Comput. Geom.}
{\bf 19}, no. 3, Special Issue, 461--469, (1998)


 
\end{thebibliography}

\end{document}